\documentclass[11pt]{article}
\usepackage{amsfonts}
\usepackage{amssymb}
\usepackage{amsmath}
\usepackage[all]{xy}

\def\Ker{\mathop{\rm Ker}\nolimits}
\def\Hom{\mathop{\rm Hom}\nolimits}
\def\Coker{\mathop{\rm Coker}\nolimits}
\def\fd{\mathop{\rm fd}\nolimits}
\def\pd{\mathop{\rm pd}\nolimits}
\def\id{\mathop{\rm id}\nolimits}
\def\rfpd{\mathop{\rm r.fpd}\nolimits}

\def\lFID{\mathop{\rm l.FID}\nolimits}
\def\lFWD{\mathop{\rm l.FWD}\nolimits}
\def\rfpd{\mathop{\rm r.fpd}\nolimits}
\def\Ext{\mathop{\rm Ext}\nolimits}
\def\Tr{\mathop{\rm Tr}\nolimits}
\def\mod{\mathop{\rm mod}\nolimits}
\def\Mod{\mathop{\rm Mod}\nolimits}
\title{\Large \bf Syzygy Modules and Injective Cogenerators for Noether Rings
\thanks{2010 Mathematics Subject Classification: 16E05, 16E10,
16E65, 16P40.}
\thanks{Keywords: syzygy modules, Auslander-type condition, injective cogenerator.}\thanks{This research was partially supported by NSFC (Grant Nos 11201220 11126169, 11101217) and Nanhua University Startup Foundation for Doctors(Grant No 2010XQD32)}}
\author{Chonghui Huang\thanks{{\it E-mail address}: huangch78@163.com} \\ {\small College of Mathematics and
Physics, University of South China, Hengyang 421001, China}}

\date{}
\begin{document}
\baselineskip=18pt \maketitle

\begin{abstract}
In this paper, we focus on $n$-syzygy modules and the injective cogenerator determined by the minimal injective resolution of a noether ring. We study the properties of $n$-syzygy modules and a category $R_n(\mod R)$ which includes the category consisting of all $n$-syzygy modules and their applications on Auslander-type rings. Then, we investigate the injective cogenerators determined by the minimal injective resolution of $R$. We show that $R$ is Gorenstein with finite self-injective dimension at most $n$ if and only if $\id R\leq n$ and $\fd \bigoplus_{i=0}^n I_i(R)< \infty$. Some known results can be our corollaries.
\end{abstract}

\section {Introduction}

In this paper, $R$ is always a noether ring if not mentioned. $\Mod R$ (resp.
$\Mod R^{op}$) denotes the category of all left(resp.
right) $R$-modules; $\mod R$ (resp. $\mod R^{op}$) denotes the category of finitely generated left(resp.
right) $R$-modules; For a left $R$-module $M$, $\pd M$, $\fd M$ and $\id M$ denotes the projective, flat and injective dimension of $M$, respectively; $I_i(M)$ denotes the $i+1$-th term in the minimal injective resolution of $M$ for any $i\geq 1$.

\vspace{0.2cm}

The existence of "{\bf Hilbert's Syzygy Theorem}" make syzygy modules category be an important research object in studying properties of rings and modules. The {\bf first motivation} of this paper is the works on syzygy modules and Auslander-type rings done by M. Auslander, I. Reiten, Z.Huang and O. Iyama among other ([FGR], [AR2,3], [Hu1,2], [Iy1,2], [HI],[HH]). Let's recall some definitions.

\vspace{0.2cm}

{Definition 1.1} ([AB]) (1) {\it $M\in \mod R$ is called to be an $n$-syzygy module (of $A\in \mod R$), if there exists an exact sequence $0\to M\to P_{n-1}\to \cdots\to P_1\to P_0\to A\to 0$ with any $P_i$ is projective module in $\mod R$. $M$ is an $\infty$-syzygy module if it is $n$-syzygy module for any $n\geq 1$}.

(2) {\it A module $M\in \mod R$ is called to be $n$-torsionfree (resp. $\infty$-torsionfree) if \linebreak $\Ext_{R^{op}}^i(\Tr M, R)=0$ for any $1\leq i\leq n$ (resp. $i\geq 1$), where $\Tr M$ is the Auslander transpose of $M$}.

\vspace{0.2cm}

Denoted by $\Omega_n(\mod R)$ (resp. $\Omega_{\infty}(\mod R)$) the category consisting of all finitely generated $n$-syzygy (resp. $\infty$-syzygy) modules;  Denoted by $\mathcal{T}_n(\mod R)$ (resp. $\mathcal{T}_{\infty}(\mod R)$) the category consisting of all finitely generated $n$-torsionfree (resp. $\infty$-torsionfree) modules. It is known that $\mathcal{T}_i(\mod R)\subset \Omega_i(\mod R)$ for any $i\geq 1$. These categories play important role in the following Auslander-type rings.

\vspace{0.2cm}

{Definition 1.2} ([Iy1]) For any $n, k \geq 0$, {\it $R$ is said to be $g_n(k)$ if $\Ext_{R^{op}}^j(\Ext_{R}^{i+k}(M,$ $ R), R)=0$ for any $M\in \mod R$ and $1\leq i\leq n$ and $0\leq j\leq i-1$; and $R$ is said to be $g_n(k)^{op}$ if $R^{op}$ is $g_n(k)$. $R$ is said to be $G_n(k)$ if $R$ has a minimal injective resolution with $\fd I_i(R^{op}) <i+k$ for any $0\le i\le n-1$; and $R$ is said to be $G_n(k)^{op}$ if $R^{op}$ is $G_n(k)$}.

\vspace{0.2cm}

It follows from [Iy1, 6.1] that $R$ is $g_n(k)$ (resp. $g_n(k)^{op}$) if $R$ is $G_n(k)$ (resp. $G_n(k)^{op}$). If $R$ satisfies $G_n(0)$-condition (resp. $G_\infty(0)$-condition, $R$ is called to be an $n$-Gorenstein ring (resp. {\it Auslander ring})([FGR]). These rings are important research object in homological algebra and the representation theory of algebras. Moreover, they are closely related to many famous homological conjectures, such as {\bf Nakayama Conjecture} (every artin algebra with infinite dominant dimension is self-injective.), {\bf Auslander-Reiten Conjecture} (every Auslander artin algebra is Gorenstein.), {\bf Gorenstein Symmetry Conjecture} (the left and right self-injective dimension is equal for any artin algebra.) and so on.

\vspace{0.2cm}

In section 2, we study the injective resolutions of $n$-syzygy modules. And we find that the minimal injective resolution of a syzygy module is quite similar to that of $R$, that is, let $0\to A \to I_0(A) \to I_1(A) \to I_2(A)\to \cdots $ be a minimal injective
resolution of a $n$-syzygy module $A$, then, $I_i(A)\underset{\oplus}\leq \coprod(\underset{j\leq i}\oplus I_j(R))$ for any $0\le i\le n-1$. As an application, we get a characterization of rings satisfying $G_n(k)$ conditions. That is, $R$ is $G_n(k)^{op}$ if and only if $\fd I_i(A)\leq i+k$ for any $0\le i\le n-1$ and any $A\in \Omega_n(\mod R)$. We guess that studying the properties of such $R$-modules whose minimal injective resolution is similar to that of $R$ will help us to understand the properties of $R$ better. So, we introduce the notion of modules admitting $R_n$-property as following

\vspace{0.2cm}

{Definition 1.3} {\it Let $0\to M \to I_0(M) \to I_1(M) \to I_2(M)\to \cdots $ be a minimal injective
resolution of $M$. $M$ is called to have $R_n$-property (resp. $R_\infty$-property) if $I_i(M)\underset{\oplus}\leq \coprod(\underset{j\leq i}\oplus I_j(R))$ for any $0\le i\le n-1$ (resp. $i\geq 0$).}

\vspace{0.2cm}

Denoted by $R_n(\mod R)$ (resp. $R_\infty(\mod R)$) the category consisting of all modules which have $R_n$-property (resp. $R_\infty$-property). Obviously, every $n$-syzygy module belongs to $R_n(\mod R)$. We show that $R_n(\mod R)$ is resolving, that is, includes all projective modules, extension closed and closed under kernel of epimorphism. Consequently, we get $\mathcal{T}_n(\mod R)\subset \Omega_n(\mod R)\subset R_n(\mod R)$. Auslander and Bridger showed that if $\mathcal{T}_i(\mod R)=\Omega_i(\mod R)$ for any $1\le i\le n$, then $R$ is $g_n(1)$ ([AB, Proposition 2.26]). So, it is nature to ask: What will happen if $\Omega_i(\mod R)= R_i(\mod R)$ (resp. $\mathcal{T}_i(\mod R)= R_i(\mod R)$) for any $1\le i\le n$. In Theorem 2.5, we show that $\Omega_i(\mod R)= R_i(\mod R)$ for any $1\le i\le n$ if and only if $\mathcal{T}_i(\mod R)= R_i(\mod R)$ for any $1\le i\le n$, in this case, $R$ is $G_n(1)^{op}$. This result partly answers the above question. And there are examples showing that the converse does not hold true. Then, it is nice to study when will $\Omega_i(\mod R)= R_i(\mod R)$ hold for any $1\le i\le n$. On the other hand, as is known to all, let $A$ be a finitely generated $R$-module, it is difficult to judge whether $A$ is an $n$-syzygy module. For an artin algebra $R$, we show that if $R$ is $G_n(0)$, then $\Omega_i(\mod R)= R_i(\mod R)$ for any $1\le i\leq n$ (Theorem 2.6). This result give us a new way to judge whether a module is $n$-syzygy module when $R$ is a $k$-Gorenstein artin algebra for any $n\leq k$.

\vspace{0.2cm}

On the other hand, J-I. Miyachi has a series of result about the injective cogenerators determined by the minimal injective resolution of a noether ring $R$ ([Mi, Theorems 2.2, 2.3, 3.9]). As applications, we give a new proof to Miyachi's results about cogenerator for category consisting of finitely projective dimension at most $n$ (Proposition 2.7). Then we show that $\bigoplus_{i=0}^n I_i(R)$ is an injective cogenerator for category consisting of finitely flat dimension at most $n$. This result generalizes [Mi, Theorem 3.9].

\vspace{0.2cm}

{\bf Another motivation} of this paper is the works on injective cogenerators done by Y. Iwanaga, J-I. Miyachi among others. It is known to all that the minimal injective resolutions of rings (resp. modules) determine many properties of them. Among those lots of properties of the minimal injective resolutions of rings, it is interesting to study the injective cogenerators of $\mod R$ (resp. $\Mod R$). A famous question concerned with the the minimal injective resolutions of rings and the cogenerators for $\mod R$ still remains open: Is $\underset{i\geq 0}\oplus I_i$ a cogenerator for $\mod R$ over an artin algebra $R$?

\vspace{0.2cm}

In [Iw1, Theorem 2], Iwanaga shows that if $\id R^{op}\leq n$, then $\underset{i\leq n}\oplus I_i(R)$ is an injective cogenerator for $\Mod R$. In [HH1, Lemma 4.5], Huangs give a wider version: If $^{\bot _n}{_RR}$ has the torsionless property, $\bigoplus_{i=0}^n I_i(R)$ is an injective cogenerator for $\Mod R$. But all these results didn't give the sufficient condition of $\underset{i\leq n}\oplus I_i(R)$ be an injective cogenarator for $\Mod R$.

\vspace{0.2cm}

The {\bf last aim} of this paper is to continue the work of Z. Huang and C. Huang in [HH1]. Denoted by $^{\bot}{_RR}=\{M\in \mod R\ |\ \Ext
_R^i(_RM, {_RR})=0$ for any $i \geq 1\}$ (resp. $^{\bot}R_R=\{N\in \mod R^{op}\ |\ \Ext^i_{R^{op}}(N_R, {R_R})=0$ for any $i \geq 1\}$). $M$ is said to have {\it Gorenstein dimension zero} if $M\in {^{\bot}{_RR}}$ and $\Tr M \in {^{\bot}R_R}$. Denoted by $\mathcal {G}(\mod R)$ the category consisting of all finitely generated modules which have Gorenstein dimension zero. It is known that $\mathcal {G}(\mod R)\subset \mathcal{T}_{\infty}(\mod R)\subset \Omega_{\infty}(\mod R)\subset R_{\infty}(\mod R)$. Huangs proved in [HH1, Theorem 1.4] that $R$ is Gorenstein with self-injective dimension at most $n$ if and only if every finitely generated module has finite Gorenstein dimension at most $n$, if and only if $\mathcal{G}(\mod R)\supset \Omega_{n}(\mod R)$. Moreover, they also proved that $\id R^{op}\leq n$ if every finitely generated module has an $\infty$-torsionfree resolution with finite lengthen at most $n$. Then, it is nature to ask: What is $R$ if every finitely generated module has an $\infty$-syzygy (resp. $R_{\infty}(\mod R)$) resolution with finite lengthen at most $n$?

\vspace{0.2cm}

In section 3, we give a necessary and sufficient condition for $\underset{i\leq n}\oplus I_i(R)$ being an injective cogenarator for $\Mod R$ (Theorem 3.1). We prove that: $\underset{i\leq n}\oplus I_i(R)$ is an injective cogenarator for $\Mod R$ if and only if $R_n(\mod R)=R_{n+1}(\mod R)=R_{\infty}(\mod R)$. As a corollary, we have that $\id R\leq n$ yields $\underset{i\leq n}\oplus I_i(R)$ is an injective cogenarator for $\Mod R$ and $\underset{i\leq n}\oplus I_i(R^{op})$ is an injective cogenarator for $\Mod R^{op}$. This result generalizes the one of Iwanaga ([Iw1, Theorem 2]). Consequently, we proved that: $\id R\leq n$ and $\fd \bigoplus_{i=0}^n I_i(R)< \infty$ if and only if $R$ is Gorenstein with self-injective dimension at most $n$ (Theorem 3.3). Finally, we give a partial answer to the above question and show that: If every finitely generated module has an $\infty$-syzygy resolution with finite length at most $n$, then $\underset{i\leq n}\oplus I_i(R)$ is an injective cogenerator of $\Mod R$ and $\lFID R\leq n$ (hence $\rfpd R\leq n$).

{\section {Syzygy modules over noether ring}}

\vspace{0.2cm}

Firstly, we give a result which plays role in this section.

\vspace{0.2cm}

{\bf Lemma 2.1} {\it Let $0\to A\to B\to C\to 0$ be an exact sequence in $\mod R$ and put $I_{-1}(M)=0$ for any $M\in \mod R$. Then  }

(1){\it $I_i(B)\underset{\oplus}\leq I_i(A)\oplus I_i(C)$ for any $i\geq 0$. }

(2) {\it $I_i(C)\underset{\oplus}\leq I_i(B)\oplus I_{i+1}(A)\oplus (I_{i-1}(A)\oplus I_{i-1}(C))$ for any $i\geq 0$.}

(3) {\it $I_i(A)\underset{\oplus}\leq I_i(B)\oplus (I_{i-1}(A)\oplus I_{i-1}(C))$ for any $i\geq 0$.}

{Proof}. (1) is trivial by horseshoe lemma.

(2) We first proof the case $i=0$. We have the following push-out diagram: $$\xymatrix{
   &  & 0 \ar[d]
   & 0 \ar[d]  &  \\
  0 \ar[r] & A \ar@{=}[d]\ar[r] & B \ar[d]
  \ar[r] & C \ar[d] \ar[r] & 0  \\
  0  \ar[r] & A \ar[r] & I_0(B) \ar[d]
  \ar[r] & D \ar[d] \ar[r] & 0\\
  &  & \Omega_{-1}(B) \ar[d]
  \ar@{=}[r]& \Omega_{-1}(B) \ar[d] &  \\
 &  & 0  & 0  &   }$$ and $$\xymatrix{
   &  & 0 \ar[d]
   & 0 \ar[d]  &  \\
  0 \ar[r] & A \ar@{=}[d]\ar[r] & I_0(A) \ar[d]
  \ar[r] & \Omega_{-1}(A) \ar[d] \ar[r] & 0  \\
  0  \ar[r] & A \ar[r] & I(B) \ar[d]
  \ar[r] & D \ar[d] \ar[r] & 0\\
  &  & I'_0 \ar[d]
  \ar@{=}[r]& I'_0 \ar[d] &  \\
 &  & 0  & 0  &   }$$ where $D=\Coker(A\to I_0(B))$ and $I'_0=\Coker(I_0(A)\to I_0(B))$. It is easy to see $I'_0$ is a direct summand of $I_0(B)$. The fact $C\subset D$ yields that $I_0(C)\underset{\oplus}\leq I_0(D)$. By the exactness of the sequence $0\to \Omega_{-1}(A)\to D\to I'_0\to 0$, we get $I_0(D) \underset{\oplus}\leq I_1(A)\oplus I'_0$. Then by (1), we have $I_0(C)\underset{\oplus}\leq I_0(B)\oplus I_{1}(A)$.

Now suppose $i\geq 1$. Using horseshoe lemma and put $D_1=\Coker(B\to I_0(A)\oplus I_0(C))$. We have the following push-out diagram: $$\xymatrix{
   &  & 0 \ar[d]
   & 0 \ar[d]  &  \\
  0 \ar[r] & B \ar@{=}[d]\ar[r] & I_0(B) \ar[d]
  \ar[r] & \Omega_{-1}(B) \ar[d] \ar[r] & 0  \\
  0  \ar[r] & B \ar[r] & I_0(A)\oplus I_0(C) \ar[d]
  \ar[r] & D_1 \ar[d] \ar[r] & 0\\
  &  & I_0 \ar[d]
  \ar@{=}[r]& I_0 \ar[d] &  \\
 &  & 0  & 0  &   }$$where $I_0=\Coker(I_0(B)\to I_0(A)\oplus I_0(C))$. It is clear that $I_0$ is a direct summand of $I_0(A)\oplus I_0(C)$. Then by (1),
 we have $I_0(D_1)\underset{\oplus}\leq I_1(B)\oplus I_0(A)\oplus I_0(C)$. Moreover, the exactness of $0\to \Omega_{-1}(B)\to D_1\to I_0\to 0$ and horseshoe lemma yield that $I_i(D_1)\underset{\oplus}\leq I_{i+1}(B)$ for any $i\geq 1$. Now we consider the following exact sequence:$$0\to \Omega_{-1}(A)\to D_1\to \Omega_{-1}(C)\to 0$$
Then by induction, $I_i(C)=I_{i-1}(\Omega_{-1}(C))\underset{\oplus}\leq I_{i}(\Omega_{-1}(A))\oplus I_{i-1}(D_1)\oplus (I_{i-2}(\Omega_{-1}(A))\oplus I_{i-2}(\Omega_{-1}(C)))$. Hence we get $I_i(C)\underset{\oplus}\leq I_{i+1}(A)\oplus I_i(B)\oplus (I_{i-1}(A)\oplus I_{i-1}(C))$.

(3) The case $i=0$ is trivial. The exact sequence $0\to \Omega_{-1}(A)\to D_1\to \Omega_{-1}(C)\to 0$ yields $I_1(A)\underset{\oplus}\leq I_{0}(D_1)\underset{\oplus}\leq I_1(B)\oplus I_0(A)\oplus I_0(C)$. Replace the exact sequence $0\to A\to B\to C\to 0$ by exact sequence $0\to \Omega_{-1}(A)\to D_1\to \Omega_{-1}(C)\to 0$, we have $I_2(A)\underset{\oplus}\leq I_{1}(\Omega_{-1}(A))\underset{\oplus}\leq I_1(D_1)\oplus I_1(A)\oplus I_1(C)$. Repeat this program and applying the fact $I_i(D_1)\underset{\oplus}\leq I_{i+1}(B)$, we get $I_i(A)\underset{\oplus}\leq I_i(B)\oplus (I_{i-1}(A)\oplus I_{i-1}(C))$ for any $i\geq 1$. So we have done. \hfill $\square$

\vspace{0.2cm}

As a direct application of Lemma 2.1, we consider the minimal injective resolution of any $n-$syzygy module and get the following result.

\vspace{0.2cm}

{\bf Proposition 2.2} {\it Let $A$ be any $n-$syzygy module in $\mod R$. Then $I_i(A)$ is a direct summand of the direct sum of $\underset{j\leq i}\oplus I_j(R)$ for any $0\le i\leq n-1$. }

\vspace{0.2cm}

{Proof}. Let $M$ be any module in $\mod R$ and $\cdots\to P_{n-1}\to P_{n-2}\to \cdots\to P_0\to M\to 0$ be any projective resolution of $M$ in $\mod R$.
Put $A_1=\Ker(P_{0}\to M)$ and $A_i=\Ker(P_{i-1}\to P_{i-2})$ for any $i\geq 2$. Then it is clear that $I_0(A_1)\underset{\oplus}\leq \coprod I_0(R)$. By Lemma 2.1(3), we have $I_0(A_2)\underset{\oplus}\leq I_0(P_1)\underset{\oplus}\leq \coprod I_0(R)$ and $I_1(A_2)\underset{\oplus}\leq I_1(P_1)\oplus (I_{0}(A_2)\oplus I_{0}(A_1))\underset{\oplus}\leq I_1(P_1)\oplus (I_{0}(P_1)\oplus I_{0}(P_0))\underset{\oplus}\leq \coprod(\underset{j\leq 1}\oplus I_j(R))$. It is not difficult to induce that $I_i(A_n)\underset{\oplus}\leq \coprod(\underset{j\leq i}\oplus I_j(R))$ for any $i<n$.\hfill $\square$

\vspace{0.2cm}

By Proposition 2.2, $R$ is $G_n(k)^{op}$ yields that $\fd I_i(A)\leq i+k$ for any $0\le i\leq n-1$ and any $A\in \Omega_n(\mod R)$. Moreover, every $R$-module with Gorenstein dimension zero in $\mod R$ is an $i$-syzygy for any $i\geq 1$. So we have the following

\vspace{0.2cm}

{\bf Corollary 2.3} {\it $R$ is $G_n(k)^{op}$ if and only if $\fd I_i(A)\leq i+k$ for any $0\le i\le n-1$ and any $A\in \Omega_n(\mod R)$ if and only if $\fd I_i(G)\leq i+k$ for any $0\le i\le n-1$ and any $A\in \mathcal {G}(\mod R)$}.

\vspace{0.2cm}

It is natural that $R_n(\mod R)\subset R_{n+1}(\mod R)$, $R_n(\mod R)$ is closed under extension, direct sum and direct summand. Recall from [AR1] that a subcategory $\mathcal{X}$ of $\mod R$ is called to be {\it resolving} if $\mathcal{X}$ includes all finitely generated projective $R$-module, closed under extension and closed under kernel of epimorphism. These category play important part in approximation theory of modules. By Proposition 2.2, every project module in $\mod R$ belongs to $R_n(\mod R)$. Moreover, we can see that $R_n(\mod R)$ is closed under kernel of epimorphism by Lemma 2.1(3). Actually, we have a stronger result.

\vspace{0.2cm}

{\bf Proposition 2.4} {\it Let $0\to A\to B\to C\to 0$ be an exact sequence in $\mod R$. Then $B\in R_{n}(\mod R)$ and $C\in R_{n-1}(\mod R)$ yield that $A\in R_{n}(\mod R)$. Hence $R_{n}(\mod R)$ is resolving.}

\vspace{0.2cm}

By the above statements, we get $\mathcal{T}_n(\mod R)\subset \Omega_n(\mod R)\subset R_n(\mod R)$. As Auslander and Reiten showed in [AR3] that $\mathcal{T}_i(\mod R)=\Omega_i(\mod R)$ for any $1\le i\le n$ if and only if $R$ is $g_n(1)$. So it is interesting to discuss the following question: what is $R$ if $R_i(\mod R)=\Omega_i(\mod R)$ for any $1\leq i\leq n$?

\vspace{0.2cm}

By Proposition 2.4, $R_i(\mod R)$ is resolving (so it is extension closed) for any $i\geq 1$. Then by [AR3, Theorem 1.7], we see $\Omega_i(\mod R)$ is extension closed for any $1\leq i\leq n$ if and only if $\Omega_i(\mod R)$ is extension closed and $\Omega_i(\mod R)=\mathcal{T}_i(\mod R)$ for any $1\leq i\leq n$. So we have the following result which answer the above question.

\vspace{0.2cm}

{\bf Theorem 2.5} {\it $R_i(\mod R)=\Omega_i(\mod R)$ for any $1\leq i\leq n$ if and only if $R_i(\mod R)=\mathcal{T}_i(\mod R)$ for any $1\leq i\leq n$, in this case, $R$ is $G_n(1)^{op}$.}

\vspace{0.2cm}

{Remark}: The inverse statement of Theorem 2.5 is not true in generally. That is, $R$ is $G_n(1)^{op}$ can't derive the fact $R_i(\mod R)=\Omega_i(\mod R)$ for any $1\leq i\leq n$. Let $\Lambda$ be an artin algebra given by the
following quiver Q:$$\xymatrix{1 \ar[r]^{a_1} &
2\ar[r]^{a_2}\ar[d]^{a_3}&3\\
 & 4& }$$ with the relation ${a_2a_1=0=a_3a_1}$. Then $\pd_{\Lambda}I_0(\Lambda)=1$. By the Definition 1.3, $I_0(\Lambda)\in R_1(\mod \Lambda)$. But it don't belong to $\Omega_i(\mod \Lambda)$ (or else, $I_0(\Lambda)$ is projective-injective). The facts $\Lambda\in \Omega_i(\mod \Lambda)$ and $R_i(\mod \Lambda)\subseteq R_1(\mod \Lambda)$ for any $i\geq 1$ means that $\Omega_i(\mod \Lambda)\subseteq \mathrm R_i(\mod \Lambda)$ does not hold true in generally. Then it is nature to ask: When will $R_i(\mod R)=\Omega_i(\mod R)$ for any $1\leq i\leq n$ hold true?

\vspace{0.2cm}

 Recall that the left orthogonal class of a subcategory $\mathcal C$, denoted by $^\perp \mathcal C$, is defined to be the subcategory of $\mod R$ consisting of those $X$ such that $\Ext^i_R(X, {\mathcal C})$ for any $i\geq 1$. The right orthogonal class of a subcategory $\mathcal D$, denoted by $ \mathcal D^\perp$, is defined to be the subcategory of $\mod R$ consisting of those $Y$ such that $\Ext^i_R({\mathcal C}, Y)$ for any $i\geq 1$. For any subcategory $\mathcal C$ and $\mathcal D$, we have $^\perp \mathcal C$ is resolving and $ \mathcal D^\perp$ is coresolving.

\vspace{0.2cm}

Let $\mathcal X$ be a full subcategory of $\mod R$. Recall that a morphism $f:\ X\to C$ is said to be a right $\mathcal X$-approximation if $\Hom_R(X', X)\to \Hom_R(X', C)\to 0$ is exact for any $X'\in \mathcal X$. $\mathcal X$ is said to be {\it contravariantly finite} in $\mod R$ if every $C$ in $\mod R$ admits a right $\mathcal X$-approximation. Denoted by ${inj}_n$ the subcategory consisting of all module with injective dimension at most $n$, then it is quit easy to check that $(\Omega_i(\mod R))^\perp=\ { {inj}_i} $ for any $i\geq 1$. We have the following

\vspace{0.2cm}

{\bf Theorem 2.6} {\it Let $R$ be an $n$-Gorenstein artin algebra. Then $R_i(\mod R)=\Omega_i(\mod R)$ for any $1\leq i\leq n$.}

\vspace{0.2cm}

{Proof}. Firstly, we claim that $R_n(\mod R)^\perp\supseteq {inj}_n$ for any $n\geq 1$. Let $N$ be any module in ${inj}_n$ and $A$ be any module in $R_i(\mod R)$. Since $R$ is $n$-Gorenstein, we have $\pd I_i(R)\leq i$ for any $0\leq i\leq n-1$. So it is nature that $\Ext_R^{j+i}(I_i(A), N)=0$ for any $j\geq 1$ and $i\geq 0$. Therefore, we get the following $$\Ext_R^{j}(A, N)\cong\Ext_R^{j+1}(\Omega_{-1}(A), N)\cong\cdots \cong\Ext_R^{j+n}(\Omega_{-(n+1)}(A), N)\cong 0\ {\rm for\ any\ } j\geq 1.$$ So we get $R_n(\mod R)^\perp\supseteq {inj}_n$ for any $n\geq 1$.

Suppose $i\geq 1$. By Proposition 2.2, we have $\Omega_i\subseteq R_i(\mod R)$. So we get $R_i(\mod R)^\perp \subseteq \Omega_i^\perp= {inj}_i$. But we also have $R_i(\mod R)^\perp \supseteq {inj}_i$. Then we have $R_i(\mod R)^\perp={inj}_i$ and hence $^\perp(R_i(\mod R)^\perp)=\ ^\perp({inj}_i)$.

If $R$ is $n$-Gorenstein, then we can see $\Omega_i(\mod R)$ is contravariantly finite and resolving for any $1\leq i\leq n$. Then by [AR1, Proposition 3.3], we have $\Omega_i(\mod R)=\ {^\perp {inj}_i}$ for any $1\leq i\leq n$. Hence we have $R_i(\mod R)
\subset\ ^\perp(R_i(\mod R)^\perp)\subset \Omega_i(\mod R)$ for any $1\leq i\leq n$. That is, $R_i(\mod R)=\Omega_i(\mod R)$ for any $1\leq i\leq n$. \hfill$\square$

\vspace{0.2cm}

{Remark}: Let $R$ be an $n$-Gorenstein artin algebra. It is known that syzygy modules play important role in characterizing properties of these rings. Another fact is that it is quite difficult to judge whether a module $A$ is a $k$-syzygy module for some $k$. [AR3] shows that if $R$ is $g_n(1)$, then $\Omega_i(\mod R)=\mathcal{T}_i(\mod R)$ for any $1\leq i\leq n$. But this method is not directly enough and it has less operability. This theorem actually give a new way to judge whether a module $A$ is a $i$-syzygy module for any $1\leq i\leq n$ and it is more convenient to study the minimal injective resolution and projective resolution.

\vspace{0.2cm}

Let $\mathcal A$ be an full subcategory of $\Mod R$. Recall that an object $X\in \Mod R$ is called a $\sum$-embedding cogenerator (resp. a finitely embedding cogenerator, a cogenerator) for $\mathcal A$ provided that every object in $\mathcal A$ admits an injection to some direct sum (resp. finite direct sum, direct product) of copies of $X$ in $\mathcal A$. Miyachi proved the following two results by derived categories. (1) Let $R$ be a left coherent ring. Then $\bigoplus_{i=0}^n I_i(R)$ is a finitely embedding cogenerator for the category of finitely present left $R$-modules of projective dimension at most $n$. (2) Let $R$ be a left noether ring. Then $\bigoplus_{i=0}^n I_i(R)$ is a $\sum$-embedding cogenerator for the category of left $R$-modules of projective dimension at most $n$. As application of Lemma 2.1(2), we give a directly proof to the his results as following.

\vspace{0.2cm}

{\bf Proposition 2.7} ([Mi, Proposition 2.1 and Theorem 2.2]) {\it (1) Let $R$ be a left coherent ring. Then $\bigoplus_{i=0}^n I_i(R)$ is a finitely embedding cogenerator for the category of finitely present left $R$-modules of projective dimension at most $n$.}

\vspace{0.2cm}

{\it (2) Let $R$ be a left noether ring. Then $\bigoplus_{i=0}^n I_i(R)$ is a $\sum$-embedding cogenerator for the category of left $R$-modules of projective dimension at most $n$.}

\vspace{0.2cm}

{Proof}. (1) Because $R$ is a left coherent ring, every $A\in \mod R$ admits a finitely generated projective resolution. Let $0\to P_n\to P_{n-1}\to \cdots\to P_0\to A\to 0$ be a projective resolution in $\mod R$ for $A\in \mod R$ with $\pd A\leq n$ and any direct sum of injective modules is still injective. Then, by Lemma 2.1(2) and Proposition 2.2, we have the following
\begin{eqnarray*}
\nonumber I_0(A)&\underset{\oplus}\leq &I_0(P_0) \oplus I_1(\Omega_1(A)) \\
\nonumber   &\underset{\oplus}\leq&I_0(P_0)\oplus [I_1(P_1)\oplus I_2(\Omega_2(A))\oplus (I_0(\Omega(A_1))\oplus I_0(P_1))]\\
\nonumber & \underset{\oplus}\leq  & [\coprod (I_0(R)\oplus I_1(R))]\oplus I_2(\Omega_2(A)) \\
\nonumber &\cdots \cdots& \\
\nonumber & \underset{\oplus}\leq & [\coprod (I_0(R)\oplus I_1(R)\oplus\cdots \oplus I_{n-1}(R))]\oplus I_n(\Omega_n(A))
\end{eqnarray*} Since we have $\Omega_n(A)=P_n$, we get $I_0(A)\underset{\oplus}\leq  \coprod (I_0(R)\oplus I_1(R)\oplus\cdots \oplus I_{n}(R))$. Then we get the assertion by the facts that $A\in \mod R$ and $P_i\mod R$ for any $0\leq i\leq n$.

(2) The proof is an analog to  that of (1), so we omit it. \hfill{$\square$}

\vspace{0.2cm}

Moreover, let $k$ be a commutative ring and $R$ a left noether projective $k$-algebra. Miyachi also proved that $\bigoplus_{i=0}^n I_i(R)$ is an injective cogenerator for the category of left $R$-modules of flat dimension at most $n$ ([Mi, Theorem 3.9]). We get a stronger version of his result as following

\vspace{0.2cm}

{\bf Theorem 2.8} {\it Let $R$ be a left noether ring. Then $\bigoplus_{i=0}^n I_i(R)$ is an injective cogenerator for the category of left $R$-modules of flat dimension at most $n$.}

\vspace{0.2cm}

{Proof}. Let $F$ be any flat left $R$-module and $A$ be any left $R$-module with flat dimension at most $n$. Then $M$ is a direct limit of a direct system of finitely generated projective right $R$-modules $\{ Q_j\}_{j\in J}$. By the proof of [Hu3, Theorem 4.1], it is known that $I_k(F)\cong \underset{j\in J}{\underset{\rightarrow}{\rm lim}}I_k(Q_j)$. Because every $I_k(Q_j)$ is a direct summand of a direct sum  of $I_k(R)$, denote it by $\coprod_F I_k(R)$   applies the same techniques to that of $I_k(F)\cong \underset{j\in J}{\underset{\rightarrow}{\rm lim}}I_k(Q_j)$, we have $I_k(F) \underset{\oplus}\leq\underset{j\in J}{\underset{\rightarrow}{\rm lim}}\coprod_F I_k(R)=\coprod_F I_k(R)$ for any $k\geq 0$, where $I$ is a direct index set. Then by a similar statement to the proof of Proposition 2.7(1), we get
$I_0(A)\underset{\oplus}\leq\coprod \bigoplus_{i=0}^n I_i(R)$. That is, $\bigoplus_{i=0}^n I_i(R)$ is an injective cogenerator for $\Mod R$.\hfill{$\square$}

\vspace{0.5cm}

\section {Injective cogenerators determined by minimal injective resolution of ring}

\vspace{0.2cm}

In this section, we throw our sight on injective cogenerators determined by minimal injective resolution of $R$ and resolutions of modules.

\vspace{0.2cm}

{\bf Theorem 3.1} {\it $\bigoplus_{i=0}^n I_i(R)$ is an injective cogenerator for
$\Mod R$ if and only if $R_n(\mod R)=R_{n+1}(\mod R)=R_{\infty}(\mod R)$.}

\vspace{0.2cm}

{Proof.} The sufficiency is trivial. Now suppose $R_n(\mod R)=R_{n+1}(\mod R)$. By Proposition 2.2 and the assumption, $\Omega_n(A)\in R_{n+1}(\mod R)$ for any $A\in \mod R$. Using Proposition 2.1(2) repeatedly, we get $I_0(S)\underset{\oplus}\leq \coprod(\underset{i\leq n}\oplus I_i(R))$ for any simple left $R$-module $S$. Hence $\Hom _{R}(S, \bigoplus_{i=0}^n I_i(R))\neq 0$. Thus we conclude that $\bigoplus_{i=0}^n I_i(R)$ is
an injective cogenerator for $\Mod R$ and $R_n(\mod R)=R_{\infty}(\mod R)$.\hfill{$\square$}

\vspace{0.2cm}

Notice that the condition $\id R\leq n$ yields that $R_n(\mod R)=R_{\infty}(\mod R)$. So the following corollary actually generalize [Iw1, Theorem 2].

\vspace{0.2cm}

{\bf Corollary 3.2} {\it If $\id R\leq n$, then $\bigoplus_{i=0}^n I_i(R)$ is an injective cogenerator for
$\Mod R$ and $\bigoplus_{i=0}^n I_i(R^{op})$ is an injective cogenerator for $\Mod R^{op}$}.

\vspace{0.2cm}

As far as the {\bf GSC} is concerned, we have the following

\vspace{0.2cm}

{\bf Theorem 3.3} {\it The following conditions are equivalent.}

(1) {\it $\id R\leq n$ and $\fd \bigoplus_{i=0}^n I_i(R)< \infty$.}

(2) {\it $R$ is Gorenstein with self-injective dimension at most $n$. }

\vspace{0.2cm}

{Proof}. We need only to show $(1)\Rightarrow (2)$. By Corollary 3.2, we get that $\bigoplus_{i=0}^n I_i(R)$ is an injective cogenerator for
$\Mod R$. So $\Mod R$ has an injective cogenerator with finite flat dimension by the assumption. Hence we have $\id ^{op}R< \infty$ by [Iw2, Proposition 1]. Then by [Za, Lemma A], we have $\id R^{op}=\id R\leq n$. We have done.  \hfill{$\square$}

\vspace{0.2cm}

{\bf Corollary 3.4} {\it For any $n\geq 1$ and $l\geq 0$,  $\id R\leq n$ and $R$ is $G_n(l)^{op}$ $\Leftrightarrow$ $R$ is Gorenstein with self-injective dimension at most $n$.}

\vspace{0.2cm}

Especially, recall that Auslander and Reiten proved that if $R$ is an artin Auslander algebra (i.e. $R$ is $n$-Gorenstein for any $n\geq 1$), then $\id R\leq \infty$ yields that $R$ has already been Gorenstein [AR2, Corollary 5.5]. Note that the condition '$R$ is $n$-Gorenstein' is left-right symmetric. We improve their result as following

\vspace{0.2cm}

{\bf Corollary 3.5} {\it If $R$ is an $n$-Gorenstein ring with $\id R\leq n$, then $R$ is already Gorenstein.}

\vspace{0.2cm}

It is known from [HH1, Theorem 1.4] that if every $n$-syzygy module in $\mod R$ has Gorenstein dimension zero, then $R$ is already a Gorenstein ring with self-injective dimension at most $n$. Note that $\mathcal{T}_n(\mod R)\subset \Omega_n(\mod R)$, then by [HH1, Proposition 4.1], we have the following

{\bf Proposition 3.6} {\it The following statement are equivalent.}

(1) {\it Every $n$-torsionfree module in $\mod R$ has Gorenstein dimension zero.}

(2) {\it Every $n$-torsionfree module in $\mod R^{op}$ has Gorenstein dimension zero.}

(3) {\it Both ${^{\bot_n}{_RR}}$ and ${^{\bot_n}{R_R}}$ have the torsionless property (that is, ${^{\bot_n}{_RR}}\subseteq \mathcal{T}_1(\mod R)$ and ${^{\bot_n}{R_R}}\subset \mathcal{T}_1(\mod R^{op})$).}

\vspace{0.2cm}

Assume that $n$ and $k$ are
positive integers and $^{\bot_n}{_RR}$ has the torsionless property. Huangs also proved that if $R$ is $g_n(k)$ or $g_n(k)^{op}$ , then $\id _{R^{op}} R \leq n+k-1$ (see [HH1, Theorem 4.7]). If every $n$-torsionfree module in $\mod R$ (resp. $\mod R^{op}$) has Gorenstein dimension zero and $R$ is $g_n(k)$ or $g_n(k)^{op}$, by Proposition 3.6 and [HH1, Theorem 4.7], we have $\id R\leq 2n-1$ and $\id R^{op}\leq n+k-1$, that is, we have the following

\vspace{0.2cm}

{\bf Proposition 3.7.} {\it Assume that $n$ and $k$ are positive integers. If every $n$-torsionfree module in $\mod R$ (resp. $\mod R^{op}$) has Gorenstein dimension zero and $R$ is $g_n(k)$ or $g_n(k)^{op}$, then $R$ is Gorenstein with self-injective dimension at most $n+k-1$}

\vspace{0.2cm}

It is known that $R$ is Gorenstein with self-injective dimension at most $n$ if and only if every finitely generated $R$-module (resp. $R^{op}$-module) has finite Gorenstein dimension at most $n$; and $R^{op}$ has finite self-injective dimension at most $n$ if every finitely generated $R$-module has finite torsionfree dimension at most $n$ (see [HH1, Theorems 1.4 and 3.6]). Moreover, It is known from Theorem 3.1 that if every finitely generated $R$-module has a $R_\infty$ resolution with finite length at most $n$ if and only if $\bigoplus_{i=0}^n I_i(R)$ is an injective cogenerator for $\Mod R$. Notice that it is known to all that $\mathcal {G}(\mod R)\subset \mathcal{T}_{\infty}(\mod R)\subset \Omega_{\infty}(\mod R)\subset R_{\infty}(\mod R)$. So, it interesting to study the $\infty$-syzygy dimension of modules. Now we introduce the notion of the $\infty$-syzygy dimension of modules as follows.

\vspace{0.2cm}

{Definition 3.8.} For a module $M\in \mod R$, the {\it
$\infty$-syzygy dimension} of $M$, denoted by $\Omega_{\infty}-\dim _RM$, is
defined as inf$\{ n\geq 0\ |\ $there exists an exact sequence $0 \to
X_{n} \to \cdots \to X_{1} \to X_{0} \to M \to 0$ in $\mod R$ with
all $X_i\in \Omega(\mod R)\}$. We set $\Omega_{\infty}-\dim _RM$
infinity if no such integer exists.

\vspace{0.2cm}

By Proposition 2.2, it is nature that $\Omega_{\infty}(\mod R)\subset R_{\infty}(\mod R)$. So, by Proposition 2.4 and Theorem 3.1, we have the following

\vspace{0.2cm}

{\bf Corollary 3.9.} {\it Let $n$ be a non-negative integer. Then every finitely generated $R$-module has finite $\infty$-syzygy dimension at most $n$ yields that $\bigoplus_{i=0}^n I_i(R)$ is an injective cogenerator for $\Mod R$.}

\vspace{0.2cm}

Actually, we have a more strong answer. Recall that the finitistic dimension of a ring $R$ are defined as follows:

$\lFID R=sup \{\id M\ |\ M\in\Mod R\ with\ \id M< \infty\ \}$

$\lFWD R=sup \{\fd M\ |\ M\in\Mod R\ with\ \fd M< \infty\ \}$

$\rfpd R=sup \{\pd M\ |\ M\in\mod R^{op}\ with\ \pd M< \infty\ \}$

\vspace{0.2cm}

It is known that $\rfpd R\leq \lFWD R=\lFID R$. Moreover, it is not difficult to see that if $\rfpd R\leq n$, then $\bigoplus_{i=0}^n I_i(R)$ is an injective cogenerator for $\Mod R$. We gave a strong version of Corollary 3.9 as following

\vspace{0.2cm}

{\bf Proposition 3.10.} {\it Let $R$ be any Noether ring and $n$ be a negative integer. If every finitely generated $R$-module has finite $\infty$-syzygy dimension at most $n$, then $\rfpd R\leq n$ and hence $\bigoplus_{i=0}^n I_i(R)$ is an injective cogenerator for $\Mod R$.}

\vspace{0.2cm}

{Proof}. Let $0\to X_n\to \cdots\to X_0\to A\to 0$ be a $\infty$-syzygy resolution of $A\in\mod R$ and $M$ be any $R$-module with $\id M=m<\infty$. Because $\Omega_{\infty}(\mod R)\subset \Omega_m(\mod R)$ and $\Omega_m(\mod R)^\bot=inj_m$, we get $\Ext^i_R(X_i,M)=0$ for any $i\geq 1$ and $j\geq 0$. So we have $\Ext^i_R(A,M)=0$ for any $i\geq n+1$. That is, $\lFID R\leq n$ and hence $\rfpd R\leq n$. \hfill{$\square$}

\end{document}